\documentclass{amsart}

\usepackage{amssymb,amsmath}

\newtheorem{theorem}{Theorem}
\newtheorem{lemma}[theorem]{Lemma}

\newtheorem{proposition}[theorem]{Proposition}

\theoremstyle{definition}
\newtheorem{definition}[theorem]{Definition}
\newtheorem{example}[theorem]{Example}

\theoremstyle{remark}
\newtheorem{remark}[theorem]{Remark}

\newcounter{th}

\newcommand{\C}{{\mathcal C}}

\newcommand{\cP}{{\mathcal P}}

\newcommand{\bC}{{\mathbb{C}}}
\newcommand{\bN}{{\mathbb{N}}}

\newcommand{\ov}[1]{\overline{#1}}

\newcommand{\Tr}{\mathrm{Tr}}

\newcommand{\jed}{{\mathbb{I}}}

\title[Positive maps between $M_2(\bC)$ and $M_n(\bC)$]{On 
the structure of positive maps \\ between matrix algebras}
\author[W. A. Majewski]{W{\l}adys{\l}aw A. Majewski}
\address{Institute of Theoretical Physics and Astrophysics, Gda{\'n}sk University,
Wita Stwosza 57, 80-952 Gda{\'n}sk, Po\-land}
\email{fizwam@univ.gda.pl}
\author[M. Marciniak]{Marcin Marciniak}
\address{Institute of Theoretical Physics and Astrophysics, Gda{\'n}sk University,
Wita Stwosza 57, 80-952 Gda{\'n}sk, Po\-land}
\email{matmm@univ.gda.pl}
\keywords{Positive maps, decomposable maps, face structure}
\subjclass[2000]{47B65, 47L07}
\thanks{W.A.M. is supported by EU grant SCALA FP6-2004-IST no. 015714 while M.M. is supported 
by the MNiSW research grant P03A 013 30.}
\begin{document}
%
%
%
%

\begin{abstract}
A partial description of the structure of positive unital maps 
$\phi: M_2(\bC) \to M_{n+1}(\bC)$ ($n\geq 2$) is given.
\end{abstract}

\maketitle

\section{Introduction}
We will be concerned with linear positive maps $\phi:M_m(\bC)\to M_n(\bC)$. 
To fix notation we begin with setting up the notation and the relevant terminology (cf. \cite{MM3}).
We say that $\phi$ is {\it positive} if $\phi(A)$ is a positive element in
$M_n(\bC)$ for every positive matrix from $M_m(\bC)$. If $k\in\bN$,
then $\phi$ is said to be {\it $k$-positive} (respectively
{\it $k$-copositive}) whenever $[\phi(A_{ij})]_{i,j=1}^k$ (respectively
$[\phi(A_{ji})]_{i,j=1}^k$) is positive in $M_k(M_n(\bC))$ for
every positive element $[A_{ij}]_{i,j=1}^k$ of $M_k(M_m(\bC))$. If
$\phi$ is $k$-positive (respectively $k$-copositive) for every
$k\in\bN$ then we say that $\phi$ is {\it completely positive}
(respectively {\it completely copositive}\/). Finally, we say that the map
$\phi$ is {\it decomposable} if it has the form
$\phi=\phi_1+\phi_2$ where $\phi_1$ is a completely
positive map while $\phi_2$ is a completely copositive one.

By $\cP(m,n)$ we denote the set of all positive maps acting between
$M_m(\bC)$ and $M_n(\bC)$ and by $\cP_1(m,n)$ -- the subset of
$\cP(m,n)$ composed of all positive unital maps (i.e. such that
$\phi(\jed)=\jed$). Recall that $\cP(m,n)$ has the structure of a
convex cone while $\cP_1(m,n)$ is its convex subset.

In the sequel we will use the notion of a face of a convex cone.
\begin{definition}\label{face}
Let $C$ be a convex cone. We say that a convex subcone $F\subset C$
is a \textit{face of $C$} if for every $c_1,c_2\in C$ the condition
$c_1+c_2\in F$ implies $c_1,c_2\in F$.

A face $F$ is said to be a \textit{maximal face}\/ if $F$ is a proper subcone
of $C$ and for every face $G$ such that $F\subseteq G$ we have $G=F$
or $G=C$.
\end{definition}

The following theorem of Kye gives a nice characterization of
maximal faces in the cone $\cP(m,n)$.
\begin{theorem}[\cite{Kye}]\label{maxface}
A convex subset $F\subset\cP(m,n)$ is a maximal face of $\cP(m,n)$
if and only if there are vectors $\xi\in\bC^m$ and $\eta\in\bC^n$
such that $F=F_{\xi,\eta}$ where
\begin{equation}\label{F}
F_{\xi,\eta}=\{\phi\in\cP(m,n):\,\phi(P_\xi)\eta=0\}
\end{equation}
and $P_\xi$ denotes the one-dimensional orthogonal projection in
$M_m(\bC)$ onto the subspace generated by the vector $\xi$.
\end{theorem}

The aim of this paper is to go one step further in clarification of the structure
of positive maps between $M_2(\bC)$ and $M_n(\bC)$. 
It is worth pointing out that many open problems
in quantum computing demand the better knowledge of this structure.
Consequently, our results shed new light on the structure of positive 
maps as well as on the nature of entanglement
(cf. \cite{MM1}, and for relation to quantum correlations see \cite{Maj}).

We recall (see \cite{S,W}) that every elements of
$\cP(2,2)$, $\cP(2,3)$ and $\cP(3,2)$ are decomposable. 
Contrary, $\cP(n,m)$ with $m,n\geq 3$ contains nondecomposable maps.
In \cite{MM2} we proved that if $\phi$ is extremal element
of $\cP_1(2,2)$ then its decomposition is unique. Moreover, we provided
a full description of this decomposition. 
In the case $m>2$ or $n>2$ the problem of finding decomposition is still
unsolved. In this paper we consider the next step for partial solution of this 
problem, namely for the case $m=2$ and $n\geq 3$.
Our approach will be based on the method
of the so called Choi matrix. 

To give a brief exposition of this method, we
recall (see \cite{Choi1,MM1} for
details) that if $\phi:M_m(\bC)\to M_n(\bC)$ is a linear map and
$\{E_{ij}\}_{i,j=1}^m$ is a system of matrix units in $M_m(\bC)$,
then the matrix
\begin{equation}\label{Cm}
\mathbf{H}_\phi=[\phi(E_{ij})]_{i,j=1}^m\in
M_m(M_n(\bC)),\end{equation} is called the Choi matrix of $\phi$
with respect to the system $\{E_{ij}\}$. Complete positivity of
$\phi$ is equivalent to positivity of $\mathbf{H}_\phi$ while
positivity of $\phi$ is equivalent to block-positivity of
$\mathbf{H}_\phi$ (see [Choi1], [MM1]). A matrix $[A_{ij}]_{i,j=1}^m\in M_m(M_n(\bC))$
(where $A_{ij}\in M_n(\bC)$) is called {\it block-positive} if 
$\sum_{i,j=1}^m\overline{\lambda_i}\lambda_j\langle\xi,A_{ij}\xi\rangle\geq 0$ for any $\xi\in\bC^n$
and $\lambda_1,\ldots,\lambda_m\in\bC$.

It was shown in Lemma 2.3 in \cite{MM2} that the general form of the Choi matrix
of a positive map $\phi$ belonging to some maximal face of $\cP(2,2)$ is the following
\begin{equation}
\label{Choi22}
\mathbf{H}_\phi=\left[\begin{array}{cc|cc}a&c&0&y\\\overline{c}&b&\overline{z}&t\\\hline
0&z&0&0\\\overline{y}&\overline{t}&0&u\end{array}\right].
\end{equation}
Here $a,b,u\geq 0$, $c,y,z,t\in\bC$ and the following inequalities are satisfied:
\begin{enumerate}
\item[(I)] $|c|^2\leq ab$,
\item[(II)] $|t|^2\leq bu$,
\item[(III)] $|y|+|z|\leq (au)^{1/2}$.
\end{enumerate}
It will turn out that in the case $\phi:M_2(\bC)\to M_{n+1}(\bC)$, $n\geq 2$,
the Choi matrix has the form which is
similar to (\ref{Choi22}) but some of the coefficients have to be matrices (see \cite{MM3}).
The main results of our paper is an analysis of the Tang's maps in the Choi matrix setting
and proving some partial results about the
structure of positive maps in the case $\phi:M_2(\bC)\to M_{n+1}(\bC)$.

\section{$\cP(2,n+1)$ maps and Tang's maps}
In this section we summarize without proofs the relevant material 
on the Choi matrix method for $\cP(2,n+1)$ (see \cite{MM3}) and we indicate how this technique may be used
to analysis of nondecomposable maps.
Let $\{e_1,e_2\}$ and $\{f_1,f_2,\ldots,f_{n+1}\}$
denote the standard orthonormal bases of the spaces $\bC^2$ and
$\bC^{n+1}$ respectively, and let $\{E_{ij}\}_{i,j=1}^2$ and
$\{F_{kl}\}_{k,l=1}^{n+1}$ be systems of matrix matrix units in
$M_2(\bC)$ and $M_{n+1}(\bC)$ associated with these bases. We assume
that $\phi\in F_{\xi,\eta}$ for some $\xi\in\bC^2$ and
$\eta\in\bC^{n+1}$. By taking the map $A\mapsto V^*\phi(WAW^*)V$
for suitable $W\in U(2)$ and $V\in U(n+1)$ we can assume without
loss of generality that $\xi=e_2$ and $\eta=f_1$. Then the Choi
matrix of $\phi$ has the form
\begin{equation}\label{Cm1}
\mathbf{H}=\left[\begin{array}{cccc|cccc}
a&c_1&\ldots&c_n&x&y_1&\ldots&y_n \\
\ov{c_1}&b_{11}&\ldots&b_{1n}&\ov{z_1}&t_{11}&\ldots &t_{1n}\\
\vdots &\vdots & &\vdots &\vdots &\vdots & &\vdots \\
\ov{c_n}&b_{n1}&\ldots &b_{nn}&\ov{z_n}&t_{n1}&\ldots &t_{nn}\\\hline
\ov{x}&z_1&\ldots &z_n&0&0&\ldots&0\\
\ov{y_1}&\ov{t_{11}}&\ldots&\ov{t_{n1}}&0&u_{11}&\ldots&u_{1n}\\
\vdots&\vdots& &\vdots&\vdots&\vdots& &\vdots\\
\ov{y_n}&\ov{t_{1n}}&\ldots&\ov{t_{nn}}&0&u_{n1}&\ldots&u_{nn}
\end{array}\right]
\end{equation}
We introduce the following notations:
$$C=\left[\begin{array}{ccc}c_1&\ldots&c_n\end{array}\right],\quad
Y=\left[\begin{array}{ccc}y_1&\ldots&y_n\end{array}\right],\quad
Z=\left[\begin{array}{ccc}z_1&\ldots&z_n\end{array}\right],$$
$$B=\left[\begin{array}{ccc}b_{11}&\ldots&b_{1n}\\\vdots& &\vdots\\b_{n1}&\ldots&b_{nn}\end{array}\right],\quad
T=\left[\begin{array}{ccc}t_{11}&\ldots&t_{1n}\\\vdots&
&\vdots\\t_{n1}&\ldots&t_{nn}\end{array}\right],\quad
U=\left[\begin{array}{ccc}u_{11}&\ldots&u_{1n}\\\vdots&
&\vdots\\u_{n1}&\ldots&u_{nn}\end{array}\right].$$ 
The matrix
(\ref{Cm1}) can be rewritten in the following form
\begin{equation}\label{Cm2}
\mathbf{H}=\left[\begin{array}{cc|cc} a&C&x&Y \\ C^*&B&Z^*&T
\\\hline \ov{x}&Z&0&0 \\ Y^*&T^*&0&U\end{array}\right].
\end{equation}
The symbol $0$ in the right-bottom block has three different
meanings. It denotes
$0$, $\left[\begin{array}{ccc}0&\ldots&0\end{array}\right]$ or
$\left[\begin{array}{c}0\\\vdots\\0\end{array}\right]$
respectively. We have the following
\begin{proposition}[\cite{MM3}]\label{Cmpos}
Let $\phi:M_2(\bC)\to M_{n+1}(\bC)$ be a positive map with the
Choi matrix of the form (\ref{Cm2}). Then the following relations
hold:
\begin{enumerate}
\item $a\geq 0$ and $B$, $U$ are positive matrices,
\item if $a=0$ then $C=0$, and if $a>0$ then $C^*C\leq aB$,
\item $x=0$,
\item the matrix
$\left[\begin{array}{c|c}B&T\\\hline T^*&U\end{array}\right]\in
M_2(M_n(\bC))$ is block-positive. 
\end{enumerate}
\end{proposition}
In the example below, we will be concerned with the two-parameter family
of nondecomposable maps (cf. \cite{Tang}). Here the important point to note is 
the fact that $\cP(2,4)$ and $\cP(3,3)$
are the lowest dimensional cases having nondecomposable maps.
Therefore the detailed analysis of such maps should yield necessary informations for
explanations of the occurrence of nondecomposability.
\begin{example}
Let $\phi_0:M_2(\bC)\to M_4(\bC)$ be a linear map defined by
\begin{equation}\label{phiT}
\phi_0\left(\left[\begin{array}{cc}a&b\\c&d\end{array}\right]\right)=
\left[\begin{array}{cccc}
(1-\varepsilon)a+\mu^2d & -b & \mu c & -\mu d \\
-c & a+2d & -2b & 0 \\
\mu b & -2c & 2a+2d & -2b \\
-\mu d & 0 & -2c & a+d
\end{array}\right],
\end{equation}
where $0<\mu<1$ and $0<\varepsilon\leq\frac{1}{6}\mu^2$.
It is proved in \cite{Tang} that $\phi_0$ is nondecomposable.
One can see that $\phi_0$ has the following Choi matrix 
\begin{equation}\label{ChoiT}
H_{\phi_0}=\left[\begin{array}{cccc|cccc}
1-\varepsilon & 0 & 0 & 0 & 0 & -1 & 0 & 0 \\
0 & 1 & 0 & 0 & 0 & 0 & -2 & 0 \\
0 & 0 & 2 & 0 & \mu & 0 & 0 & -2 \\
0 & 0 & 0 & 1 & 0 & 0 & 0 & 0 \\ \hline
0 & 0 & \mu & 0 & \mu^2 & 0 & 0 & -\mu \\
-1 & 0 & 0 & 0 & 0 & 2 & 0 & 0 \\
0 & -2 & 0 & 0 & 0 & 0 & 2 & 0 \\
0 & 0 & -2 & 0 & -\mu & 0 & 0 & 1
\end{array}\right].
\end{equation}
Observe that 
\begin{equation*}
\phi_0(\jed)=\left[\begin{array}{cccc}
1-\varepsilon+\mu^2 & 0 & 0 & -\mu \\
0 & 3 & 0 & 0 \\
0 & 0 & 4 & 0 \\
-\mu & 0 & 0 & 2
\end{array}\right].
\end{equation*}
Let $\rho=\sqrt{1-\varepsilon+\mu^2}$ and 
\begin{equation*}
\delta=\left|\begin{array}{cc}
1-\varepsilon+\mu^2 & -\mu \\
-\mu & 2
\end{array}\right|^{1/2}=
\sqrt{2-2\varepsilon+\mu^2}.
\end{equation*}
Then
$\phi_0(\jed)^{-1/2}$ is of the form
\begin{equation*}
\phi_0(\jed)^{-1/2}=
\left[\begin{array}{cccc}
\dfrac{\beta}{\delta} & 0 & 0 & -\dfrac{\gamma}{\delta} \\
0 & \dfrac{1}{\sqrt{3}} & 0 & 0 \\
0 & 0 & \dfrac{1}{2} & 0 \\
-\dfrac{\gamma}{\delta} & 0 & 0 & \dfrac{\alpha}{\delta}
\end{array}\right]
\end{equation*}
where 
and $\alpha,\beta >0$, $\gamma\in\mathbb{R}$ are such that
\begin{eqnarray}
\alpha^2+\gamma^2 & = & \rho^2 \nonumber\\
\beta^2+\gamma^2 & = & 2 \label{abc}\\
(\alpha+\beta)\gamma & = & -\mu.\nonumber
\end{eqnarray}

Let us define $\phi_1:M_2(\bC)\to M_4(\bC)$ by 
$$\phi_1(A)=\phi_0(\jed)^{-1/2}\phi_0(A)\phi_0(\jed)^{-1/2},\quad A\in M_2(\bC).$$
Then
\begin{eqnarray*}
\phi_1(E_{11}) & = & 
\left[\begin{array}{cccc}
\dfrac{(1-\varepsilon)\beta^2+\gamma^2}{\delta^2} & 0 & 0 & -\dfrac{[(1-\varepsilon)\beta+\alpha]\gamma}{\delta^2} \\
0 & \dfrac{1}{3} & 0 & 0 \\
0 & 0 & \dfrac{1}{2} & 0 \\
-\dfrac{[(1-\varepsilon)\beta+\alpha]\gamma}{\delta^2} & 0 & 0 & \dfrac{(1-\varepsilon)\gamma^2+\alpha^2}{\delta^2}
\end{array}\right]  \\[6mm]
\phi_1(E_{22}) & = & 
\left[\begin{array}{cccc}
\dfrac{(\mu \beta+\gamma)^2}{\delta^2} & 0 & 0 & -\dfrac{(\mu \beta+\gamma)(\mu \gamma+\alpha)}{\delta^2} \\
0 & \dfrac{2}{3} & 0 & 0 \\
0 & 0 & \dfrac{1}{2} & 0 \\
-\dfrac{(\mu \beta+\gamma)(\mu \gamma+\alpha)}{\delta^2} & 0 & 0 & \dfrac{(\mu \gamma+\alpha)^2}{\delta^2}
\end{array}\right] \\[6mm]
\phi_1(E_{12})& = &
\left[\begin{array}{cccc}
0 & -\dfrac{\beta}{\delta\sqrt{3}} & 0 & 0 \\
0 & 0 & -\dfrac{1}{\sqrt{3}} & 0 \\
\dfrac{\mu \beta+2\gamma}{2\delta} & 0 & 0 & -\dfrac{\mu \gamma+2\alpha}{2\delta} \\
0 & \dfrac{\gamma}{\delta\sqrt{3}} & 0 & 0 
\end{array}\right]
\end{eqnarray*}

One can deduce from (\ref{abc}) that
\begin{equation}\label{w}
(\mu \gamma+\alpha)^2+(\mu \beta+\gamma)^2=\rho^2.
\end{equation}
Let 
\begin{equation*}
W=\left[\begin{array}{cccc}
\dfrac{\mu \gamma+\alpha}{\sqrt{1-\varepsilon+\mu^2}} & 0 & 0 & \dfrac{\mu \beta+\gamma}{\sqrt{1-\varepsilon+\mu^2}} \\
0 & 1 & 0 & 0\\
0 & 0 & 1 & 0\\
\dfrac{\mu \beta+\gamma}{\sqrt{1-\varepsilon+\mu^2}} & 0 & 0 & -\dfrac{\mu \gamma+\alpha}{\sqrt{1-\varepsilon+\mu^2}} 
\end{array}\right].
\end{equation*}
It follows from (\ref{w}) that $W$ is a unitary matrix.
Define $\phi:M_2(\bC)\to M_4(\bC)$ by 
$\phi(A)=W^*\phi_1(A)W$. Then the Choi matrix of $\phi$ is of the form
\begin{equation}\label{Phitil}
\mathbf{H}_\phi=\left[\begin{array}{cccc|cccc}
1 & 0 & 0 & 0 & 0 & -\dfrac{1}{\sqrt{3}\rho} & 0 & 0\\
0 & \dfrac{1}{3} & 0 & 0 & 0 & 0 & -\dfrac{1}{\sqrt{3}} & 0 \\
0 & 0 & \dfrac{1}{2} & 0 & -\dfrac{\mu}{2\rho} & 0 & 0 & \dfrac{\delta}{2\rho} \\
0 & 0 & 0 & \dfrac{1-\varepsilon}{\delta^2} & 0 & -\dfrac{\mu}{\sqrt{3}\delta\rho} & 0 & 0 \\
\hline
0 & 0& -\dfrac{\mu}{2\rho} & 0 & 0 & 0 & 0 & 0 \\
-\dfrac{1}{\sqrt{3}\rho} & 0& 0 & -\dfrac{\mu}{\sqrt{3}\delta\rho} & 0 & \dfrac{2}{3} & 0 & 0 \\
0 & -\dfrac{1}{\sqrt{3}} & 0 & 0 & 0 & 0 & \dfrac{1}{2} & 0 \\
0 & 0 & \dfrac{\delta}{2\rho} & 0 & 0 & 0 & 0 & \dfrac{\rho^2}{\delta^2}
\end{array}\right]
\end{equation}
One can see that $\psi\left(\left[\begin{array}{cc}0&0\\0&1\end{array}\right]\right)
\left[\begin{array}{c}1\\0\\0\\0\end{array}\right]=0$,
so $\psi\in F_{\eta,\xi}$ (cf. Theorem 2), where
$\eta=\left[\begin{array}{c}0\\1\end{array}\right]$ and $\xi=\left[\begin{array}{c}1\\0\\0\\0\end{array}\right]$.

Observe that blocks which form the Choi matrix (\ref{Phitil}) as in (\ref{Cm2}) are of the form
\begin{eqnarray*}
a&=&1,\qquad
C=0,\qquad
Y=\left[\begin{array}{ccc}-\dfrac{1}{\sqrt{3}\delta} & 0 & 0\end{array}\right],\qquad
Z=\left[\begin{array}{ccc}0& -\dfrac{\mu}{2\rho} & 0\end{array}\right],\\
B&=&\left[\begin{array}{ccc} \dfrac{1}{3} & 0 & 0  \\
0 & \dfrac{1}{2} & 0 \\
0 & 0 & \dfrac{1-\varepsilon}{\delta^2} \end{array}\right],\quad
U=\left[\begin{array}{ccc}
\dfrac{2}{3} & 0 & 0 \\
0 & \dfrac{1}{2} & 0 \\
0 & 0 & \dfrac{\rho^2}{\delta^2}\end{array}\right],\\
T&=&\left[\begin{array}{ccc}
0 & -\dfrac{1}{\sqrt{3}} & 0 \\
0 & 0 & \dfrac{\delta}{2\rho} \\
-\dfrac{\mu}{\sqrt{3}\delta\rho} & 0 & 0 
\end{array}\right].
\end{eqnarray*}
It is worth observing that the vectors $C,Y,Z$ are orthogonal, matrices $B,U$ are diagonal,
while $T$ is ``complementary" to diagonal matrices $B$ and $U$.
This observation is useful in understanding the peculiarity of nondecomposable maps.
\end{example}

In the sequel we will need some technicalities. 
For $X=\left[\begin{array}{ccc}x_1&\ldots&x_n\end{array}\right]\in M_{1,n}(\bC)$ 
we define $\Vert X\Vert=\left(\sum_{i=1}^n|x_i|^2\right)^{1/2}$.
By $|X|$ we denote the square ($n\times n$)-matrix $(X^*X)^{1/2}$. We identify
elements of $M_{n,1}(\bC)$ with vectors from $\bC^n$ and for any $X\in M_{1,n}(\bC)$ define 
a unit vector $\xi_X\in\bC^n$ by $\xi_X=\Vert X\Vert^{-1}X^*$.
\begin{proposition}
Let $X,X_1,X_2\in M_{1,n}(\bC)$. Then 
\begin{enumerate}
\item[(1)] $|X|=\Vert X\Vert P_{\xi_X}$,
where $P_\xi$ denotes the orthogonal 
projection onto the one-di\-men\-sio\-nal subspace in $\bC^n$ generated by
a vector $\xi\in\bC^n$;
\item[(2)] $|X_1||X_2|=\langle\xi_{X_1},\xi_{X_2}\rangle X_1^*X_2$.
\end{enumerate}
\end{proposition}
\begin{proof}
(1) Let $\eta\in\bC^n$. Since $\eta$ is considered also as an element of $M_{n,1}(\bC)$ the
multiplication of matrices $X\eta$ has sense. As a result we obtain a $1\times 1$-matrix
which can be interpreted as a number. With this identifications we have the equality
$$X\eta=\langle X^*,\eta\rangle$$
where $X^*$ on the right hand side is considered as a vector from $\bC^n$, and $\langle\cdot,\cdot\rangle$
denotes the usual scalar product in $\bC^n$. Now we can calculate 
$$\langle\eta,X^*X\eta\rangle=\langle X\eta,X\eta\rangle=\Vert X\eta\Vert^2=|\langle X^*,\eta\rangle|^2=
\Vert X\Vert^2|\langle \xi_X,\eta\rangle|^2$$

(2) If $X_1=0$ or $X_2=0$ then the equality is obvious. In the case both $X_1$ and $X_2$ are nonzero the
equality follows from the following computations
\begin{eqnarray*}
|X_1||X_2|&=&\Vert X_1\Vert^{-1}\Vert X_2\Vert^{-1}|X_1|^2|X_2|^2=\Vert X_1\Vert^{-1}\Vert X_2\Vert^{-1}
X_1^*X_1X_2^*X_2 
\\
&=&\Vert X_1\Vert^{-1}\Vert X_2\Vert^{-1}X_1^*(X_1X_2^*)X_2=
\Vert X_1\Vert^{-1}\Vert X_2\Vert^{-1}\langle X_1^*,X_2^*\rangle X_1^*X_2 
\\
&=&
\langle\xi_{X_1},\xi_{X_2}\rangle X_1^*X_2 
\end{eqnarray*}
\end{proof}

To proceed with the study of Tang's maps we recall some general
properties of maps in $\cP(2,n+1)$ (cf. \cite{MM3}). We start with
\begin{proposition}[\cite{MM3}]\label{lemmain}
A map $\phi$ with the Choi matrix of the form
\begin{equation}\label{Cm3}
\mathbf{H}=\left[\begin{array}{cc|cc} a&C&0&Y \\ C^*&B&Z^*&T
\\\hline 0&Z&0&0 \\ Y^*&T^*&0&U\end{array}\right].
\end{equation}
is positive if and only if the inequality
\begin{equation}\label{nier1}
\left|\langle Y^*,\mathit{\Gamma}^\tau\rangle+
\ov{\langle Z^*,\mathit{\Gamma}^\tau\rangle} +
\Tr\left(\mathit{\Lambda}^\tau
T\right)\right|^2\leq
\left[\alpha a+\Tr\left(\mathit{\Lambda}^\tau
B\right)+2\Re\left\{\langle C^*,\mathit{\Gamma}^\tau
\rangle\right\}\right]\Tr\left(\mathit{\Lambda}^\tau U\right)
\end{equation}
holds for every $\alpha\in\bC$, matrices
$\mathit{\Gamma}=\left[\begin{array}{ccc}\gamma_1&\ldots&\gamma_n\end{array}\right]$
and
$\mathit{\Lambda}=\left[\begin{array}{ccc}\lambda_{11}&\ldots&\lambda_{1n}\\
\vdots&
&\vdots\\\lambda_{n1}&\ldots&\lambda_{nn}\end{array}\right]$,
$\gamma_i\in\bC$, $\lambda_{ij}\in\bC$ for $i,j=1,2,\ldots,n$, such
that
\begin{enumerate}
\item $\alpha\geq 0$ and $\mathit{\Lambda}\geq 0$,
\item $\mathit{\Gamma}^*\mathit{\Gamma}\leq\alpha\mathit{\Lambda}$.
\end{enumerate}
The superscript $\tau$ denotes the transposition of matrices.
\end{proposition}
and
\begin{theorem}[\cite{MM3}]\label{Cmyz}
If the assumptions of Proposition 3 
are fulfilled, then
\begin{equation}\label{yz}
|Y|+|Z|\leq a^{1/2}U^{1/2}. 
\end{equation}
\end{theorem}
\begin{remark}
One can easily check that the nondecomposable maps described in Example 4 fulfill the above inequality.
It is easy to check that in this case the inequality is proper (i.e. there is no equality). This
observation will be crucial for next section.
\end{remark}

As we mentioned, for $\cP(2,n)$, $n>3$, there are nondecomposable maps. The proposition below
provides the characterization of completely positive and completely copositive components
of $\cP(2,n)$.
\begin{proposition}[\cite{MM3}]\label{CP}
Let $\phi:M_2(\bC)\to M_{n+1}(\bC)$ be a linear map with the Choi
matrix of the form (\ref{Cm3}). Then
the map $\phi$ is completely positive (resp. completely copositive) if and only if the
following conditions hold:
\begin{enumerate}
\item[(1)] $Z=0$ (resp. $Y=0$),
\item[(2)] the matrix
$\left[\begin{array}{ccc}a&C&Y\\C^*&B&T\\Y^*&T^*&U\end{array}\right]$
(resp. $\left[\begin{array}{ccc}a&C&Z\\C^*&B&T^*\\Z^*&T&U\end{array}\right]$)
is a positive element of the algebra $M_{2n+1}(\bC)$.
\end{enumerate}
In particular, the condition (2) implies:
\begin{enumerate}
\item[(3)] if $B$ is an invertible matrix, then $T^*B^{-1}T\leq U$ (resp. $TB^{-1}T^*\leq U$),
\item[(4)] $C^*C\leq aB$,
\item[(5)] $Y^*Y\leq aU$ (resp. $Z^*Z\leq aU$).
\end{enumerate}
\end{proposition}
This Proposition yields information about possible splitting of a decomposable map
into completely positive and completely copositive components.
To go one step further let us make the following observation.
Let $\phi:M_m(\bC)\to M_n(\bC)$ be a decomposable map and 
$\phi=\phi_1+\phi_2$ for some completely positive $\phi_1$ and completely
copositive $\phi_2$. Then from Kadison inequality we easily obtain
\begin{equation}\label{Kad}
\phi(E_{ij})^*\phi(E_{ij})\leq\Vert\phi(\jed)\Vert\left(\phi_1(E_{ii})+\phi_2(E_{jj})\right)
\end{equation}
for $i,j=1,2,\ldots,m$.

Assume now that $\phi:M_2(\bC)\to M_{n+1}(\bC)$ has the Choi matrix of the form (\ref{Cm2}).
It follows from Proposition 9 that Choi matrices of $\phi_1$ and $\phi_2$ are respectively 
\begin{equation}\label{Cmdec}
\mathbf{H_1}=\left[\begin{array}{cc|cc} a_1&C_1&0&Y \\ C_1^*&B_1&0&T_1
\\\hline 0&0&0&0 \\ Y^*&T_1^*&0&U_1\end{array}\right],\quad
\mathbf{H_2}=\left[\begin{array}{cc|cc} a_2&C_2&0&0 \\ C_2^*&B_2&Z^*&T_2
\\\hline 0&Z&0&0 \\ 0&T_2^*&0&U_2\end{array}\right].
\end{equation}

Clearly, $\mathbf{H_1} + \mathbf{H_2} = \mathbf{H}$, where $\mathbf{H}$ is the Choi matrix 
corresponding to $\phi$. The inequality (\ref{Kad}) leads to additional relations between 
components of the Choi matrices
\begin{equation*}
\left[\begin{array}{cc}\Vert Z\Vert^2&ZT\\T^*Z^*&|Y|^2+T^*T\end{array}\right]
\leq
\Vert\phi(\jed)\Vert
\left[\begin{array}{cc}a_1&C_1\\C_1^*&B_1+U_2\end{array}\right]
\end{equation*}
and
\begin{equation*}
\left[\begin{array}{cc}\Vert Y\Vert^2&YT^*\\TY^*&|Z|^2+TT^*\end{array}\right]
\leq
\Vert\phi(\jed)\Vert
\left[\begin{array}{cc}a_2&C_2\\C_2^*&B_2+U_1\end{array}\right].
\end{equation*}
It is worth pointing out that the above inequalities give a partial answer to Choi question (cf. \cite{Choi2}).
Furthermore, turning to Tang's maps one can observe that the matrix
corresponding to $\phi(E_{ij})^*\phi(E_{ij})$ is relatively large what spoils
a possibility of decomposition of these maps.

\section{On the structure of elements of $\cP(2,n+1)$.}
Giving a full description of the situation in $\cP(2,2)$ in
\cite{MM2} we proved that if $\phi:M_2(\bC)\to M_2(\bC)$ is from a large class of 
extremal positive unital maps, then the constituent maps $\phi_1$ and $\phi_2$ are uniquely 
determined (cf. Theorem 2.7 in \cite{MM2}). We recall that the Choi matrix of such
extremal map $\phi:M_2(\bC)\to M_2(\bC)$ is of the form (cf. (\ref{Choi22}))
\begin{equation}\label{Choi22ext}
\mathbf{H}_\phi=\left[\begin{array}{cc|cc}1&0&0&y\\0&1-u&\overline{z}&t\\\hline
0&z&0&0\\\overline{y}&\overline{t}&0&u\end{array}\right],
\end{equation}
where, in particular, the following equality is satisfied (cf. (III) from Section 1)
\begin{equation}\label{eqext}
|y|+|z|=u^{1/2}.
\end{equation}

In this section, motivated by the results given in the previous section 
(we `quantized' the relations (I)-(III) given at the end of Section 1),
we consider maps $\phi:M_2(\bC)\to M_{n+1}(\bC)$.
If such a map is positive unital and $\phi\in F_{e_2,f_1}$ then its Choi matrix has the
form
\begin{equation}\label{Cmunit}
\left[\begin{array}{cc|cc} 1&0&0&Y\\0&B&Z^*&T\\\hline
0&Z&0&0\\Y^*&T^*&0&U\end{array}\right],
\end{equation} 
where $B$ and
$U$ are positive matrices such that $B+U=1$ and conditions listed in
Propositions 3 
and 6 
are satisfied. 

Our object is to examine consequences
of property  
\begin{equation}\label{eqextn}
|Y|+|Z|=U^{1/2}
\end{equation}
which for $n\geq 1$ is a natural analog of (\ref{eqext}).


First, we remind the following technical 
\begin{lemma}\label{bp}
Let $\mathbf{A}=\left[\begin{array}{cc}P&S\\S^*&Q\end{array}\right]\in M_2(M_n(\bC))$, where $P,Q,S\in M_n(\bC)$, and
$P,Q\geq 0$.
The following are equivalent:
\begin{enumerate}
\item[(i)] $\mathbf{A}$ is block-positive;
\item[(ii)] $pP+sS+\ov{s}S^*+qQ\geq 0$ for every numbers $p,q,s$ such that $p,q\geq 0$ and $|s|^2\leq pq$;
\item[(iii)] $|\langle\eta,S\eta\rangle|^2\leq\langle\eta,P\eta\rangle\langle\eta,Q\eta\rangle$ for every
$\eta\in\bC^n$.
\end{enumerate}
\end{lemma}
\begin{proof}
(i)
$\Rightarrow $(ii). 
Let $\eta\in\bC^n$. 
It follows from the definition of block-positivity (cf. [MM2])
that the matrix 
$$
\left[\begin{array}{cc}
\langle\eta,P\eta\rangle   & \langle\eta,S\eta\rangle \\
\langle\eta,S^*\eta\rangle & \langle\eta,Q\eta\rangle
\end{array}\right]
$$
is positive. 
Hence the matrix 
$$
\left[\begin{array}{cc}
\langle\eta,pP\eta\rangle        & \langle\eta,sS\eta\rangle \\
\langle\eta,\ov{s}S^*\eta\rangle & \langle\eta,qQ\eta\rangle
\end{array}\right]
$$
being a Hadamard product of two positive matrices is positive as well.
Consequently, 
$$
\langle\eta,(pP+sS+\ov{s}S^*+qQ)\eta\rangle\geq 0.
$$
Since $\eta$ is arbitrary, (ii) is proved.

(ii)$\Rightarrow$(i). 
To prove that $\mathbf{A}$ is block-positive one should show
that for any $\eta\in\bC^n$ and $\mu_1,\mu_2\in\bC$ one has
$$|\mu_1|^2\langle\eta,P\eta\rangle+2\Re\left\{\mu_1\ov{\mu_2}
\langle\eta,S\eta\rangle\right\}+
|\mu_2|^2\langle\eta,Q\eta\rangle\geq 0.$$
Observe that 
$p=|\mu_1|^2$, $q=|\mu_2|^2$, 
$s=\mu_1\ov{\mu_2}$ fulfill $p,q\geq 0$ 
and $|s|^2=pq$.
So,
$$
|\mu_1|^2\langle\eta,P\eta\rangle+2\Re\left\{\mu_1\ov{\mu_2}
\langle\eta,S\eta\rangle\right\}+
|\mu_2|^2\langle\eta,Q\eta\rangle=
\langle\eta,(pP+sS+\ov{s}S^*+qQ)\eta\rangle\geq 0.
$$

(i)$\Leftrightarrow$(iii). 
Let $\eta\in\bC^n$. The positivity of the matrix
$\left[\begin{array}{cc}\langle\eta,P\eta\rangle &\langle\eta,S\eta\rangle\\
\langle\eta,S^*\eta\rangle &\langle\eta,Q\eta\rangle\end{array}\right]$
is equivalent to non-negativity of its determinant
$\langle\eta,P\eta\rangle\langle\eta,Q\eta\rangle-|\langle\eta,S\eta\rangle|^2$.
\end{proof} 
Here we give another (cf. Proposition 6) characterisation of positive maps in the language
of their Choi matrices
\begin{proposition}\label{lempos}
Let $\phi:M_2(\bC)\to M_{n+1}(\bC)$ be a linear unital map with the Choi matrix of the form
\begin{equation}\label{Choiunit}
\left[\begin{array}{cc|cc}
1&0&0&Y\\0&B&Z^*&T\\\hline 0&Z&0&0\\Y^*&T^*&0&U
\end{array}\right]
\end{equation}
where $B,U,T\in M_n(\bC)$, $Y,Z\in M_{1,n}(\bC)$, and $B,U\geq 0$.
Then the map $\phi$ is positive if and only if
$$pB+sT+\ov{s}T^*+qU\geq 0$$
and 
\begin{equation}\label{pos}
(\ov{s}Y^*+sZ^*)(sY+\ov{s}Z)\leq p^2B+p(sT+\ov{s}T^*)+pqU
\end{equation}
for every $p,q,s\in\bC$ such that $p,q\geq 0$ and $|s|^2\leq pq$.
\end{proposition}
\begin{proof}
It follows from the definition of the Choi matrix and from (\ref{Choiunit})
that 
$$
\phi\left(\left[\begin{array}{cc}p&s\\v&q\end{array}\right]\right)=
\left[\begin{array}{cc}p&sY+vZ\\
sZ^*+vY^*&pB+sT+vT^*+qU\end{array}\right].
$$
So, the map $\phi$ is positive if and only if the matrix
\begin{equation}\label{mat}
\left[\begin{array}{cc}p&sY+\ov{s}Z\\
sZ^*+\ov{s}Y^*&pB+sT+\ov{s}T^*+qU\end{array}\right].
\end{equation}
is a positive element of $M_{n+1}(\bC)$ for numbers $p,q,s$ such that $p,q\geq 0$ and
$|s|^2\leq pq$ (i.e. such that the matrix $\left[\begin{array}{cc}p&s\\\ov{s}&q\end{array}\right]$
is positive in $M_2(\bC)$). The positivity of the matrix (\ref{mat}) is equivalent to points (1) and (2)
from the statement of the lemma.
\end{proof} 

The following generalizes Lemma 8.10 from \cite{S}.
\begin{proposition}\label{TT}
Let $\phi$ be a positive unital map with the Choi matrix (\ref{Choiunit}).
Assume that $B$ is invertible.
Then the matrix
\begin{equation}\label{Tmat}
\left[\begin{array}{cc}
2B& T\\T^*&U-|Y|^2-|Z|^2
\end{array}\right]
\end{equation}
is block-positive.
\end{proposition}
\begin{proof}
Let $\eta\in\bC^n$, $\eta\neq 0$, and $p,q,s\in\bC$ be numbers such that $p,q\geq 0$
and $|s|^2=pq$.
Then from (\ref{pos}) we have
\begin{equation*}
|s|^2 \langle \eta , (|Y|^2+|Z|^2) \eta \rangle + 2 \Re \left\{s^2\langle\eta , Z^*Y \eta \rangle \right\} \leq
p^2 \langle \eta , B \eta \rangle + 2p \, \Re \left\{s \langle \eta , T \eta \rangle \right\}+ pq\langle\eta , U \eta\rangle .
\end{equation*}
Replace $s$ in this inequality by $is$ and obtain
\begin{equation*}
|s|^2 \langle \eta , (|Y|^2+|Z|^2) \eta \rangle - 2 \Re \left\{s^2 \langle \eta , Z^*Y \eta \rangle \right\} \leq
p^2\langle\eta , B \eta\rangle + 2p\, \Re \left\{is\langle\eta , T \eta\rangle \right\} + 
pq\langle\eta , U \eta\rangle .
\end{equation*}
Adding the above two inequalities one gets
\begin{equation}\label{T3}
|s|^2\langle\eta , (|Y|^2+|Z|^2) \eta\rangle\leq
p^2\langle\eta , B \eta\rangle + p\, \Re \left\{(1+i)s\langle\eta , T \eta \rangle \right\} + 
pq\langle\eta , U \eta\rangle .
\end{equation}
Let $pq=1$, and $s$ be such that $|s|=1$ and 
$\Re\left\{(1+i)s\langle\eta,T\eta\rangle\right\}=
-\sqrt{2}\,|\langle\eta,T\eta\rangle|$.
Then the inequality (\ref{T3}) takes the form
\begin{equation}\label{T4}
\langle\eta,(|Y|^2+|Z|^2)\eta\rangle\leq p^2\langle\eta,B\eta\rangle
-\sqrt{2}\,p\,\vert \langle \eta , T \eta\rangle\vert + \langle\eta , U \eta\rangle .
\end{equation}
Following the argument of St{\o}rmer in the proof of Lemma 8.10 in \cite{S}
we observe that the function
$f(x)=\langle\eta,B\eta\rangle\,x^2-\sqrt{2}\,|\langle\eta,T\eta\rangle|x+
\langle\eta,U\eta\rangle$
has its minimum for $x=2^{-1/2}\langle\eta,B\eta\rangle^{-1}|\langle\eta,
T\eta\rangle|$.
Hence, (\ref{T4}) leads to the inequality
$$\langle\eta,(|Y|^2+|Z|^2)\eta\rangle\leq -2^{-1}\langle\eta,B\eta
\rangle^{-1}|\langle\eta,T\eta\rangle|^2+\langle
\eta,U\eta\rangle$$
and finally
$$|\langle\eta,T\eta\rangle|^2\leq 2\langle\eta,B\eta\rangle\langle\eta,
(U-|Y|^2-|Z|^2)\eta\rangle.$$
By Lemma 10 this implies block-positivity of the matrix (\ref{Tmat}).
\end{proof} 

Our next results show that the property (\ref{eqextn}) in the case $n\geq 2$ has
rather restrictive consequences.
\begin{proposition}Let $\phi:M_2(\bC)\to M_{n+1}(\bC)$, $n\geq 2$, be a positive linear 
map with the Choi matrix of the form
(\ref{Choiunit}). Assume $|Y|+|Z|=U^{1/2}$. Then 
$Y$ and $Z$ are linearly dependent.
\end{proposition}
\begin{proof}
Assume on the contrary that $Y$ and $Z$ are linearly independent.
We will show that $\phi$ can not be positive in this case.
To this end 
let $p,q,s$ be numbers such that $p>0$, $q>0$ and $|s|^2\leq pq$ and define
$$D=p^2B+p(sT+\ov{s}T^*)+pqU-(\ov{s}Y^*+sZ^*)(sY+\ov{s}Z).$$
By Proposition 11 (cf. (\ref{pos})) it is enough to find numbers $p,q,s$ and a vector
$\xi_0\in\bC^n$ such that
$\langle\xi_0,D\xi_0\rangle<0$.

It follows from the assumption and Proposition 5 that
\begin{eqnarray*}
D&=&p^2B+p(sT+\ov{s}T^*)+pq(|Y|+|Z|)^2+\\
&&-\;|s|^2(|Y|^2+|Z|^2)-\ov{s}^2Y^*Z-s^2Z^*Y=\\
&=&p^2B+\left(pq-|s|^2\right)\left(|Y|^2+|Z|^2\right)+pq\left(|Y|\,|Z|+|Z|\,|Y|\right)\\
&&+\;p\left(sT+\ov{s}T^*\right)-\ov{s}^2Y^*Z-s^2Z^*Y=\\
&=&p^2B+\left(pq-|s|^2\right)\left(|Y|^2+|Z|^2\right)+psT+p\ov{s}T^*+\\
&&+ \;\left( pq\langle\xi_Y , \xi_Z \rangle - \ov{s}^2 \right) Y^*Z+ 
\left(pq\langle\xi_Z , \xi_Y \rangle - s^2 \right) Z^* Y.
\end{eqnarray*}
Let $\xi\in\bC^n$. Then 
\begin{eqnarray*}
\langle\xi,D\xi\rangle
&=&p^2\langle\xi,B\xi\rangle
+\left(pq-|s|^2\right)\langle\xi,\left(|Y|^2+|Z|^2\right)\xi\rangle + 2p\,\Re\left\{ s\langle\xi,T\xi\rangle\right\}\\
&&+\;2\Re\left\{\left( pq\langle\xi_Y,\xi_Z\rangle-\ov{s}^2 \right)\langle\xi,Y^*Z\xi\rangle\right\}=\\
&=&p^2\langle\xi,B\xi\rangle
+\left(pq-|s|^2\right)\langle\xi,\left(|Y|^2+|Z|^2\right)\xi\rangle + 2p\,\Re\left\{ s\langle\xi,T\xi\rangle\right\}\\
&&+\;2\Vert Y\Vert\Vert Z\Vert\,\Re\left\{\left(pq\langle\xi_Y,\xi_Z\rangle-
\ov{s}^2\right)\langle\xi,\xi_Y\rangle\langle \xi_Z,\xi\rangle\right\} .
\end{eqnarray*}
Let $\xi_0=\xi_Y+\xi_Z$ and $s=(pq)^{1/2}e^{i\theta}$ for some $\theta\in [0,2\pi)$. Then
\begin{eqnarray*}
\langle\xi_0,D\xi_0\rangle
&=&
p^2\langle\xi_0,B\xi_0\rangle+ 2p^{3/2}q^{1/2}\,\Re \left\{e^{i\theta}\langle\xi_0,T\xi_0\rangle\right\}\\
&&+\;
2pq\Vert Y\Vert \Vert Z\Vert\Re\left\{\left(
\langle\xi_Z,\xi_Y\rangle-
e^{-2i\theta}\right)\left(1+\langle\xi_Y,\xi_Z\rangle\right)^2\right\}.
\end{eqnarray*}
By the assumption $\xi_Y$ and $\xi_Z$ are linearly dependent. Moreover $\Vert\xi_Y\Vert=\Vert\xi_Z\Vert=1$. 
This implies that
$|\langle\xi_Z,\xi_Y\rangle|<1$, so $(1+\langle\xi_Z,\xi_Y\rangle)^2\neq 0$.
Now, choose $\theta$ such that 
$$\Re\left\{ e^{-2i\theta}(1+\langle\xi_Z,\xi_Y\rangle)^2\right\}=|1+\langle\xi_Z,\xi_Y\rangle|^2.$$
Then
\begin{eqnarray*}
\langle\xi_0,D\xi_0\rangle
&=&
p^2\langle\xi_0,B\xi_0\rangle+ 2p^{3/2}q^{1/2}\,\Re\left\{ e^{i\theta}\langle\xi_0,T\xi_0\rangle\right\}\\
&&+\;
2pq\Vert Y\Vert\,\Vert Z\Vert\,\left[
\Re\left\{\langle\xi_Y,\xi_Z\rangle(1+\langle\xi_Z,\xi_Y\rangle)^2\right\}-
|1+\langle\xi_Z,\xi_Y\rangle|^2\right].
\end{eqnarray*}
Observe that 
$$\Re\left\{\langle\xi_Y,\xi_Z\rangle(1+\langle\xi_Z,\xi_Y\rangle)^2\right\}<|1+\langle\xi_Z,\xi_Y\rangle|^2,$$
so it is possible to find $p$ sufficiently small and $q$ sufficiently large
so that $\langle\xi_0,D\xi_0\rangle$ is negative. This ends the proof.
\end{proof}
\begin{proposition}
Let $\phi:M_2(\bC)\to M_{n+1}(\bC)$ satisfy the assumptions of the previous
Proposition. If $Z=0$ and $\Vert Y\Vert<1$ (resp. $Y=0$ and $\Vert Z\Vert<1$) then $\phi$ is completely 
positive (resp. completely copositive).
\end{proposition}
\begin{proof}
It follows that $U=|Y|^2$. Moreover, the assumption $\Vert Y\Vert<1$ implies that
$B=1-|Y|^2$ is invertible. As we also have $U-|Y|^2-|Z|^2=0$,
by Proposition 12 the matrix 
$\left[\begin{array}{cc}2B&T\\T^*&0\end{array}\right]$ is block-positive. Hence $T=0$.
We conclude that the Choi matrix of $\phi$ has the form
\begin{equation*}
\left[\begin{array}{cc|cc}
1&0&0&Y\\0&1-|Y|^2&0&0\\ \hline 0&0&0&0\\Y^*&0&0&|Y|^2
\end{array}\right].
\end{equation*}
In order to finish the proof one should show (cf. Proposition 9) that the matrix
\begin{equation*}
\left[\begin{array}{ccc}1&0&Y\\0&1-|Y|^2&0\\Y^*&0&|Y|^2\end{array}\right]
\end{equation*}
is positive,
but this can be done by straightforward computations.

The proof in the case $Y=0$ follows in the same way. 
\end{proof}

As a consequence of the above results we get the following description of maps satisfying the ``quantized"
properties
(\ref{eqext}).
\begin{theorem}
Let $\phi:M_2(\bC)\to M_{n+1}(\bC)$ be a positive unital map with the Choi matrix of
the form (\ref{Choiunit}) where $|Y|+|Z|=U^{1/2}$. Then
\begin{enumerate}
\item[(1)] there are vectors $\xi\in\bC^2$ and $\eta_0\in\bC^{n+1}$ such that
\begin{equation}\label{pp}\phi\in\bigcap_{\eta\perp \eta_0}F_{\xi,\eta};\end{equation}
\item[(2)] $\phi$ is unitarily equivalent to a map with the Choi matrix of the form
\begin{equation}\label{postac}
\left[\begin{array}{ccc|ccc}
1&0&0&0&0&y\\
0&1&0&0&0&W^*\\
0&0&1-u&\ov{z}&V&t\\ \hline
0&0&z&0&0&0\\
0&0&V^*&0&0&0\\
\ov{y}&W&\ov{t}&0&0&u
\end{array}\right]\end{equation}
where in each block there are numbers on positions $[1,1]$, $[1,3]$, $[3,1]$ and $[3,3]$, 
one-row matrices from $M_{1,n-1}(\bC)$
on positions $[1,2]$ and $[3,2]$, one-column matrices from $M_{n-1,1}(\bC)$ on positions 
$[2,1]$ and $[2,3]$, and
square matrices from $M_{n-1}(\bC)$ on positions $[2,2]$.
Here $u=(|y|+|z|)^2$. Moreover, coefficients satisfy the inequality
\begin{equation}\label{rel}
|\langle\rho,Y_1^*\rangle|+|\langle\rho,Z_1^*\rangle|\leq u^{1/2}
\end{equation}
for any unit vector $\rho\in\bC^n$ where $Y_1,Z_1\in M_{1,n}(\bC)$ are defined as
$$Y_1=\left[\begin{array}{cc}\ov{y}&W\end{array}\right],\qquad
Z_1=\left[\begin{array}{cc}\ov{z}&V\end{array}\right].$$
\end{enumerate}
\end{theorem}
\begin{proof}
It follows from Proposition 13 that there is a unit vector $\eta_0\in\bC^n$ such that
$Y^*=\ov{y}\eta_0$ and $Z^*=\ov{z}\eta_0$ for some $y,z\in\bC$. Hence $|Y|=|y|P_{\eta_0}$, $|Z|=|z|P_{\eta_0}$, 
and
$U=(|y|+|z|)^2P_{\eta_0}$, where $P_{\eta_0}$ is the orthogonal projector onto the one-dimensional subspace
generated by the vector $\eta_0$. As 
$$\phi(P_{e_2})=\left[\begin{array}{cc}0&0\\0&U\end{array}\right]\in M_{n+1}(\bC)$$
then $\phi(P_{e_2})\eta=0$ for any $\eta$ orthogonal to $\eta_0$. So, from Theorem 2 we obtained (\ref{pp}).

By choosing a suitable basis of $\bC^{n+1}$ we may assume that $f_{n+1}=\eta_0$. Then
the Choi matrix (\ref{Choiunit}) takes the form
\begin{equation*}
\left[\begin{array}{ccccc|ccccc}
1&0&\cdots&0&0&0&0&\cdots&0&y\\ 
0&1&\cdots&0&0&0&t_{11}&\cdots&t_{1,n-1}&t_{1n}\\
\vdots&\vdots& &\vdots&\vdots&\vdots&\vdots& &\vdots&\vdots\\
0&0&\cdots&1&0&0&t_{n-1,1}&\cdots&t_{n-1,n-1}&t_{n-1,n}\\
0&0&\cdots&0&1-u&\ov{z}&t_{n1}&\cdots&t_{n,n-1}&t_{nn}\\\hline
0&0&\cdots&0&z&0&0&\cdots&0&0\\
0&\ov{t_{11}}&\cdots&\ov{t_{n-1,1}}&\ov{t_{n1}}&0&0&\cdots&0&0\\
\vdots&\vdots& &\vdots&\vdots&\vdots&\vdots& &\vdots&\vdots\\
0&\ov{t_{1,n-1}}&\cdots&\ov{t_{n-1,n-1}}&\ov{t_{n,n-1}}&0&0&\cdots&0&0\\
\ov{y}&\ov{t_{1n}}&\cdots&\ov{t_{n-1,n}}&\ov{t_{nn}}&0&0&\cdots&0&u
\end{array}\right].
\end{equation*}
Block-positivity of this matrix implies that the matrix
\begin{equation*}
\left[\begin{array}{ccc|ccc}
1&\cdots&0&t_{11}&\cdots&t_{1,n-1}\\
\vdots& &\vdots&\vdots& &\vdots\\
0&\cdots&1&t_{1,n-1}&\cdots&t_{n-1,n-1}\\\hline
\ov{t_{11}}&\cdots&\ov{t_{n-1,1}}&0&\cdots&0\\
\vdots& &\vdots&\vdots& &\vdots\\
\ov{t_{1,n-1}}&\cdots&\ov{t_{n-1,n-1}}&0&\cdots&0
\end{array}\right].
\end{equation*}
is also block positive, so $t_{ij}=0$ for $i,j=1,2,\ldots,n-1$. Thus we obtained 
that the Choi matrix has the form (\ref{postac}).

Now, for any $\rho\in\bC^n$, where 
$\rho=\left[\begin{array}{ccc}\rho_1&\ldots&\rho_n\end{array}\right]$,
define the following matrix from $M_{n+1,2}(\bC)$
\begin{equation*}
V_\rho=\left[\begin{array}{cccc}\ov{\rho_1}&\ldots&\ov{\rho_n}&0\\0&\ldots&0&1\end{array}\right].
\end{equation*}
One can easily check that $VV^*=1$, so a map $\psi_\rho:M_{n+1}(\bC)\to M_2(\bC):A\mapsto VAV^*$ 
is a unital and completely positive one. As a consequence we get that the map
$\psi_\rho\circ\phi:M_2(\bC)\to M_2(\bC)$ is positive and unital.
Moreover, by a straightforward calculations one can check that the Choi matrix of this map
has the form
\begin{equation*}
\left[\begin{array}{cc|cc}1&0&0&\langle\rho,Y_1^*\rangle\\
0&1-u&\ov{\langle\rho,Z_1^*\rangle}&t\\\hline
0&\langle\rho,Z_1^*\rangle&0&0\\
\ov{\langle\rho,Y_1^*\rangle}&\ov{t}&0&u
\end{array}\right].
\end{equation*}
The inequality (\ref{rel}) follows from (III) in Section 1.
\end{proof}

We end this paper by a remark that Theorem 15 gives a very useful tool for describing 
properties of extremal maps in $\cP(2,n+1)$ and it seems that following this line of research 
can give a possibility to construct some new examples of nondecomposable maps.
However, details will be contained in the forthcoming publications.


\end{document}